\documentclass{amsart}
\usepackage{amsthm}
\newtheorem{thm}{Theorem}[section]

\newtheorem{prop}[thm]{Proposition}
\theoremstyle{definition}
\newtheorem{defn}[thm]{Definition}
\theoremstyle{remark}

\numberwithin{equation}{section}

\begin{document}
\large
\title[Gabor Frames and Nonlinear Approximation]{Characterization of
Modulation Spaces by Nonlinear Approximation}%
\author{S. Samarah and S. Al-Sa'di}%

\address{Salti Samarah\\Jordan University of Science and Technology\\Jordan-Irbid\\Department of math. and Stat.\\}%
\email{samarah@just.edu.jo}%

\address{Sa'ud Al-Sa'di\\Jordan University of Science and Technology\\Jordan-Irbid\\Department of math. and Stat.\\}%
\email{saudtam@just.edu.jo}%


\begin{abstract}
 It is shown that the modulation spaces $M_{p}^{w}$ can be
 characterized by the approximation behavior of their elements
 using Local Fourier bases. In analogy to the Local Fourier bases, we show
 that the modulation spaces can also be characterized by the approximation behavior
of their elements using Gabor frames. We derive direct and inverse
approximation theorems that describe the best approximation by
linear combinations of N terms of a given function using its
modulates and translates.
\end{abstract}
\maketitle
\section{\textbf{Introduction}}

 One of the central problems of approximation theory is to
 characterize the set of functions which have a prescribed order
 of approximation by a given method of approximation. Results of
this type are known and easy to prove if we work in a Hilbert space
and the set where the approximation is seeking is an orthonormal
bases.

For general systems, there are sufficient conditions on a given
function in the space which guarantees certain rate of decrease
for the error of approximation. In this paper we consider the
method of nonlinear approximation in particular spaces called the
\emph{modulation spaces}. We investigate the approximation of
smooth functions by time-frequency shifts. The method of nonlinear
approximation has recently found many computational applications
 such as data compression, statistical estimation or adaptive
 schemes for partial differential or integral equations.

Nonlinear approximation is utilized in many numerical algorithms,
it occurs in several applications. In mathematics and applications
it is very important to write a function in some function space in
the form
\begin{equation*}
   f=\sum_{k\in\Lambda}\lambda_{k}g_{k}
\end{equation*}
where $\Lambda$ is an indexed set and $\{g_{k}\,:\,k\in\Lambda\}$ is
a set of functions. The case in which this set is obtained from a
single function is very interesting. One way to construct such set
is by using the Gabor frames. Here we characterize functions with a
given degree of nonlinear approximation for a given function $f$
belonging to some specific space using its modulates and translates.

Let $X$ be a Banach space, with norm\\ $\|.\|=\|.\|_{X}$. Then we
say that a subset $\mathcal{D}$ of functions from $X$ is a
\emph{dictionary} if :
\begin{equation}
for \;all\; g\in\mathcal{D},\;\; \|g\|_{X}=1, \;\;and\;\; g\in
\mathcal{D} \;\;implies \;\;-g\in \mathcal{D}
\end{equation}

Let $\Sigma_{N}(\mathcal{D})$ be the set of all functions in $X$
which can be written as a linear combination of at most $N$ elements
of $\mathcal{D}$, i.e.,
\begin{equation}
\Sigma_{N}(\mathcal{D})=\{s\in X;\; s=\sum_{k\in
F}c_{k}g_{k}\,,\;g_{k}\in\mathcal{D},\;c_{k}\in
\mathbb{C},\;card(F)\leq N\}
\end{equation}
From this definition we note that the sum of two elements from
$\Sigma_{N}$ is generally not in $\Sigma_{N}$, which means that the
space $\Sigma_{N}$ is not linear.

For any given $f$, the error associated to the best $N$-term
approximation to $f$ from $\mathcal{D}$ is given by:
\begin{equation}
\sigma_{N}(f,\mathcal{D})_{X}:=\sigma_{N}(f)_{X}=\inf_{s\in
\Sigma_{N}}\|f-s\|_{X}
\end{equation}

The method of nonlinear approximation was first used by Stechkin;
he characterized the space of all $f\in L_{2}(\mathbb{R})$ which
have absolutely convergent orthogonal expansion, see
\cite{Stechkin}. More precisely, he proved the following theorem:

\begin{thm}(\cite{DevoreTemlyakov}).\label{thm:Stechkinresult}
Let $\{\phi_{k};\;k\in\mathbb{ N}\}$ be an orthonormal basis for
$L_{2}(\mathbb{R})$. Then for $f\in L_{2}(\mathbb{R})$ we have:
\begin{equation*}
\sum_{k\in \mathbb{N}}|\langle f,\phi_{k}\rangle|<\infty
\;\;\Longleftrightarrow\;\;
\sum_{N=1}^{\infty}(N^{\frac{1}{2}}\sigma_{N}(f))\frac{1}{N}<\infty
\end{equation*}
\end{thm}


For a general discussion of the characterization problem, DeVore
and Temlyakove made a modification of Stechkin's result to
characterize other approximation spaces: For a general dictionary
$\mathcal{D}$, and for any $p>0$ define
\[A_{p}^{o}(\mathcal{D},M)=\{f\in \mathcal{H};\;f=\sum_{k\in\Lambda}c_{k}w_{k},
 w_{k}\in \mathcal{D},\;|\Lambda|<\infty \;and\;\left(\sum_{k\in \Lambda}|c_{k}|^{p}
  \right)^{\frac{1}{p}}\leq M \}\]
and we define $A_{p}(\mathcal{D},M)$ as the closure (in
$\mathcal{H}$) of $A_{p}^{o}(\mathcal{D},M)$. Furthermore, define
\[A_{p}(\mathcal{D})=\bigcup_{M>0}A_{p}(\mathcal{D},M)\]
and for $f\in A_{p}(\mathcal{D})$  the \emph{semi-norm}
$|f|_{A_{p}(\mathcal{D})}$ is the smallest $M$ such that $f\in
A_{p}(\mathcal{D},M)$.

Let $\mathcal{D}$ be given by an orthonormal basis
$\{\phi_{k};\,k\in \mathbb{Z}\}$. Then $f\in A_{p}(\mathcal{D})$
if and only if $ \sum_{k}|\langle f,\phi_{k}\rangle|^{p}<\infty$
and
\[|f|_{A_{p}(\mathcal{D})}= \left(\sum_{k\in \mathbb{Z}}|\langle f,\phi_{k}\rangle|^{p}\right)^{1/p}.\]

It means that we can characterize certain approximation orders by
the spaces $\mathcal{A}_{p}$. So, Stechkin's result in Theorem
(\ref{thm:Stechkinresult}) can be formulated as follows:
\[f\in \mathcal{A}_{1}(\mathcal{D}) \;\;if\;and\;only\;if\;\;
\sum_{N=1}^{\infty}(N^{1/2}\sigma_{N}(f,\mathcal{D}))\frac{1}{N} <
\infty.
\]
A slight modification of Stechkin's result due to DeVore and
Temlyakov is given in the following theorem.

\begin{thm}(\cite{DevoreTemlyakov}).
 If $\mathcal{D}$ is given by an orthonormal basis for $\mathcal{H}$, for
$\alpha>0$ and $p=(\alpha+\frac{1}{2})^{-1}$, we have :
\begin{equation}
f\in \mathcal{A}_{p}(\mathcal{D}) \;\;\Longleftrightarrow \;\;
\sum_{N=1}^{\infty}(N^{\alpha}\sigma_{N}(f))^{p}\frac{1}{N} <
\infty.
\end{equation}
\end{thm}

This theorem provides a characterization of functions with an
approximate order like $\mathcal{O}(N^{-\alpha})$, for
$\alpha=\frac{1}{p}-\frac{1}{2}$, using an orthonormal basis, and
as a special case, setting $p=1$ we get Stechkin's result. We are
now interested in characterization the approximation space
$\mathcal{A}_{p}(\mathcal{D})$ in the whole range of the
parameters $\alpha, p$ and a given dictionary $\mathcal{D}$.

\section{\textbf{Gabor Frames and Modulation Spaces}}
         \label{sec:Modul.spaces}

\ \ \ \ In this section we will collects the necessary information
the Gabor frames and modulation spaces. The modulation spaces occur
in the study of the concentration of a function in the
time-frequency plane. They were introduced in 1983 by H. Feichtinger
\cite{Feichtinger1} and were subsequently investigated in
\cite{Feichtinger,Feich and Groch3}. These spaces are defined by the
decay properties of the short time Fourier transform, and contain
many classical function spaces.

In addition to the basic definitions and notations of
\cite{Grochenig}, the Schwartz class and the space of tempered
distributions on $\mathbb{R}$ are denoted by
$\mathcal{S}(\mathbb{R})$ and $\mathcal{S'}(\mathbb{R})$
respectively. The translation and modulation operators are defined,
respectively, by:
     \begin{equation}\label{eq:translation and modulation operators}
       \textbf{\emph{T}}_{x}f(t)=f(t-x) \qquad and \qquad
       \textbf{\emph{M}}_{w}f(t)=e^{2\pi iwt}f(t)
     \end{equation}
     for $x,w\in \mathbb{R}$.
The \emph{Short-Time Fourier transform} (STFT) of a function $f\in
L^{2}(\mathbb{R})$ with
    respect to a function $g\in L^{2}(\mathbb{R})$ called the \emph{window function}, is defined as:
     \begin{equation}
       \mathcal{V}_{g}f(x,w)=\int_{\mathbb{R}}f(t)\bar{g}(t-x)e^{-2\pi
       iwt}dt=\langle
       f,\textbf{\emph{M}}_{w}\textbf{\emph{T}}_{x}g\rangle.
    \end{equation}

  A \emph{submultiplicative} weight function $v$ on $\mathbb{R}^{2}$ which is
  a positive, symmetric and continuous function and satisfies
    \begin{equation}
     v(z_{1}+z_{2})\leq v(z_{1})v(z_{2}), \qquad for\;all \;
     z_{1},z_{2}\in \mathbb{R}^{2}.
    \end{equation}
 An  $v$-\emph{moderate} weight function $m$ on $\mathbb{R}^{2}$ which is
 a positive, symmetric and continuous function and satisfies
    \begin{equation}
     m(z_{1}+z_{2})\leq C v(z_{1})m(z_{2}), \qquad for\;all \;
     z_{1},z_{2}\in \mathbb{R}^{2}.
    \end{equation}

For any continuous strictly positive function $m$ on $\mathbb{R}$,
the weighted $L_{p}$ space $L_{p}^{m}(\mathbb{R})$ is defined by the
norm
\[  \|f\|_{L_{p}^{m}}=\|fm\|_{L_{p}}
\]
and the mixed-norm spaces $L^{m}_{p,q}(\mathbb{R}^{2})$  consists of
all (Lebesgue) measurable functions on $\mathbb{R}^{2}$, such that
for a weight function $m$ on $\mathbb{R}^{2}$ the norm

 \begin{equation}
      \|F\|_{L^{m}_{p,q}}=\left(\int_{\mathbb{R}}\left( \int_{\mathbb{R}} |F(x,w)|^{p}
      \,m(x,w)^{p}dx\right)^{q/p}dw\right)^{1/q}
    \end{equation}
is finite.
 Throughout this paper, we will use two types of weights:

 Given a non-zero \emph{window function} $g\in L_{2}(\mathbb{R})$ and
   constants $\alpha,\beta>0$, the set of time-frequency
   shifts
   \begin{equation}
   \mathcal{G}(g,\alpha,\beta)=\{\mathbf{T}_{\alpha k}\mathbf{M}_{\beta n}g\,;\,k,n\in \mathbb{Z}\}
   \end{equation}
   is called a  \emph{Gabor frame} for $L_{2}(\mathbb{R})$ if there exists constants $A,B>0$
   (called frame bounds) such that for all $f\in L_{2}(\mathbb{R})$

 \begin{equation}
      A\|f\|_{L_{2}}^{2}\leq \sum_{k,n\in \mathbb{Z}}|\langle
      f,\mathbf{T}_{\alpha k}\mathbf{M}_{\beta n}g\rangle|^{2}\leq B\|f\|_{L_{2}}^{2}
 \end{equation}

%
%

Under stronger assumptions on the function $g$, the expansion in
equations (\ref{eq:f-expansion using Gabor(L_2)}) and
(\ref{eq:f-expansion 2 using Gabor(L_2)}) are valid not only in
$L_{2}$ but in the entire class of function spaces, namely, the
modulation spaces.
\begin{defn}(\cite{Grochenig}).\label{defn:Modulation Spaces}
Fix a non-zero window $g\in\mathcal{S}(\mathbb{R})$, a
$v$-moderate weight function $m$ on $\mathbb{R}^{2}$, and $1\leq
p,q\leq\infty$. Then the \emph{modulation space}
$M_{p,q}^{m}(\mathbb{R})$ consists of all tempered distributions
$f\in\mathcal{S}'(\mathbb{R})$ such that
   \begin{equation}\label{eq:M-p,q definition and V_g(f)}
     \|f\|_{M_{p,q}^{m}}=\|\mathcal{V}_{g}f\|_{L_{p,q}^{m}}=\left(\int_{\mathbb{R}}\left(\int_{\mathbb{R}}|
     \mathcal{V}_{g}f(x,w)|^{p}\,m(x,w)^{p}dx
     \right)^{q/p}dw\right)^{1/q}
   \end{equation}
is finite.
\end{defn}

Thus, $M_{p,q}^{m}$ is a Banach space of tempered distributions.
If $p=q$, then we write $M_{p}^{m}$ instead of $M_{p,p}^{m}$, and
if $m(z)\equiv 1$ on $\mathbb{R}^{2}$, then we write $M_{p,q}$ and
$M_{p}$ for $M_{p,q}^{1}$ and $M_{p}^{1}$, respectively. Some
examples of modulation spaces are the following:
\begin{enumerate}
  \item The Segal algebra $S_{0}(\mathbb{R})=M_{1,1}(\mathbb{R})$.
  \item $L_{2}(\mathbb{R})=M_{2,2}(\mathbb{R})$. However, $L_{p}$
  does not coincide with any modulation space when $p\neq 2$
  \cite{Feich-Groch-Walnut}.
  \item The Bessel potential space:\\ $H^{s}(\mathbb{R})=\{f\in \mathcal{S}';
         \|f\|_{H^{s}}=\left(\int_{\mathbb{R}}|\hat{f}(w)|^{2}\,(1+|w|^{2})^{2s}dw\right)^{1/2}<\infty \}$.
\end{enumerate}

\begin{thm}(\cite{Grochenig}).\label{thm:properties of M-p,q, Banach and invariant}
Let $m$ be a $v$-moderate weight. Then
   \begin{enumerate}
     \item $M_{p,q}^{m}(\mathbb{R})$ is a Banach space for $1\leq
           p,q\leq\infty$.
     \item $M_{p,q}^{m}$ is invariant under time-frequency
             shifts, and $\|\mathbf{T}_{x}\mathbf{M}_{w}f\|_{M_{p,q}^{m}}\leq
             Cv(x,w)\|f\|_{M_{p,q}^{m}}$.
     \item If $1\leq p_{1}, p_{2}, q_{1}, q_{2}\leq\infty$ and $m_{2}\leq Cm_{1}$,
   then
     \begin{equation*}
        M_{p_{1},q_{1}}^{m_{1}}\subseteq
             M_{p_{2},q_{2}}^{m_{2}},\; whenever\; p_{1}\leq p_{2},q_{1}\leq
            q_{2}.
    \end{equation*}
   Moreover, there exists a constant C such that
    \begin{equation}
       \|f\|_{M_{p_{2},q_{2}}^{m_{2}}}\leq C\|f\|_{M_{p_{1},q_{1}}^{m_{1}}}
    \end{equation}
    for all $f\in M_{p_{1},q_{1}}^{m_{1}}$.
   \end{enumerate}
\end{thm}

The appropriate window class in this setting is the Feichtinger
algebra
\[ M_{1}^{v}=\{f\in \mathcal{S}'(\mathbb{R}):\,\mathcal{V}_{f}f\in
L_{1}^{v}(\mathbb{R}^{2})\}.\] where $v$ is a submultiplicative
weight on $\mathbb{R}^{2}$ with polynomial growth.

\begin{thm}(\cite{Grochenig}).\label{thm:Modulation spaces and expansion, and equivalence relations}
Suppose $1\leq p,q\leq\infty$, $m$ is a $v$-moderate, $g,\gamma\in
M_{1}^{v}$. Suppose that $\{\mathbf{M}_{\beta n}
       \mathbf{T}_{\alpha k}g, \;k,n\in \mathbb{Z}\}$ generates a
       frame for $L_{2}(\mathbb{R})$,
  then for all $f\in M_{p,q}^{m}$ we have
  \begin{eqnarray}
    f&=&\sum_{n\in\mathbb{Z}}\sum_{k\in \mathbb{Z}}\langle f,\mathbf{M}_{\beta n}
       \mathbf{T}_{\alpha k}g\rangle\mathbf{M}_{\beta n}
       \mathbf{T}_{\alpha k}\gamma\\
     &=&\sum_{n\in\mathbb{Z}}\sum_{k\in \mathbb{Z}}\langle f,\mathbf{M}_{\beta n}
       \mathbf{T}_{\alpha k}\gamma\rangle\mathbf{M}_{\beta n}
       \mathbf{T}_{\alpha k}g
  \end{eqnarray}
  with unconditional convergence in $M_{p,q}^{m}$ if $p,q<\infty$,
  and weak-star convergence in $M_{\infty}^{1/v}$ otherwise.
   Furthermore, there are constants $A,B>0$ such that for all
   $f\in M_{p,q}^{m}$
   \begin{equation}\label{eq:Gabor frame condition A,B}
   A\|f\|_{M_{p,q}^{m}}\leq\left(\sum_{n\in \mathbb{Z}}\left(\sum_{k\in \mathbb{Z}}|\langle f,
       \mathbf{M}_{\beta n}\mathbf{T}_{\alpha k}g\rangle|^{p}\,m(\alpha k,\beta
       n)^{p}\right)^{q/p} \right)^{1/q}\leq B\|f\|_{M_{p,q}^{m}}
   \end{equation}
   And the norm equivalence:
   \begin{equation}
      A'\|f\|_{M_{p,q}^{m}}\leq\left(\sum_{n}\left(\sum_{k}|\langle f,\mathbf{M}_{\beta n}
       \mathbf{T}_{\alpha k}\gamma\rangle|^{p}\,m(\alpha k,\beta
       n)^{p}\right)^{q/p} \right)^{1/q}\leq B'\|f\|_{M_{p,q}^{m}}
  \end{equation}
  holds on $M_{p,q}^{m}$.
\end{thm}

\section{\textbf{Characterization of Modulation Spaces $M_{p}$}}\label{sec:charac.of Modul.spaces}
\ \ \ \ According to the definition of the approximation spaces
$\mathcal{A}_{p}(\mathcal{D})$, several questions arise:

\begin{enumerate}
  \item \emph{If $\mathcal{D}$ is given by the Fourier basis, what is
        the space that can be characterized by this basis ?.}
  \item \emph{If $\mathcal{D}$ is given by a Local Fourier basis, what
        are the spaces that can be characterized by this basis ?.}
  \item \emph{If $\mathcal{D}$ is given by Gabor frame, what are
        the spaces that can be characterized by this basis ?.}
\end{enumerate}

The first question was answered by Stechkin in Theorem
(\ref{thm:Stechkinresult}). For the second question we have the K.
Gr\"{o}chenig and S. Samarah result \cite{GrochSamarah}.

\begin{thm}(\cite{GrochSamarah}).\label{thm:iff basis Theorem}
 If $0<p<q\leq\infty$ and
$\alpha=\frac{1}{p}-\frac{1}{q}$, then:
  \begin{equation}
  f\in M_{p}^{w}\qquad if\; and\; only\; if \qquad
\sum_{N=1}^{\infty}\left(N^{\alpha}\sigma_{N}(f)_{M_{q}^{w}}
\right)^{p}\frac{1}{N}<\infty
\end{equation}
\end{thm}

In the proof of Theorem (\ref{thm:iff basis Theorem}) the basis
property was used in an essential way to rewrite the approximation
error in terms of a sequence space norm. For linearly dependent
sets it is not clear how much of Theorem (\ref{thm:iff basis
Theorem}) still holds. The next theorem prove the one-half of
Theorem (\ref{thm:iff basis Theorem}) under the weaker assumption
that the set $\{\mathbf{T}_{\beta m}\mathbf{M}_{\gamma n}\phi :
m,n\in \mathbb{Z}\}$ is a Banach frame for $M_{p}^{w}$.

\begin{prop}(\cite{Samarah}).\label{prop:one side Gabor-charac.}
 Let $\{\mathbf{T}_{\beta m}\mathbf{M}_{\gamma n}\phi : m,n\in \mathbb{Z}\}$
be a Banach frame for $M_{p}^{\omega}$ for all $0<p<\infty$. If\,
$0<p<q$, $\alpha=\frac{1}{p}-\frac{1}{q}$ and $f\in
M_{p}^{\omega}$, then
\begin{equation}
\sum_{N=1}^{\infty}\left(N^{\alpha}\sigma_{N}(f)_{M_{q}^{\omega}}
\right)^{p}\frac{1}{N}<\infty
\end{equation}
\end{prop}

Proposition (\ref{prop:one side Gabor-charac.}) does not give a
complete characterization of the modulation spaces $M_{p}$ using
the Gabor atoms as a dictionary. It proves only that for a
function $f$ in a modulation space $M_{p}$, the approximation
error has order $N^{-\alpha}$, but does not give the other
implication, i.e., if the approximation error for some function
has order $N^{-\alpha}$ for some $\alpha>0$, what does this tell
us about the space to which $f$ belongs. This inquiry will be
answered in our next work which will based on the classical
inequalities of Jackson and Bernstein where we applied with the
Gabor atoms in the modulation spaces.

Let $X$ be a Banach space in which approximation takes place and
assume that we can find a number $r>0$ and a second space $Y$
continuously embedded in $X$, and $X_{n}$ be the subsets of $X$ in
which approximants come from. Then

\emph{\textbf{Jackson Inequality}}: $\sigma_{n}(f)_{X}\leq
C\,n^{-r}\,|f|_{Y}$, \; for all $f\in Y$,\;$n=1,2,.....$

\emph{\textbf{Bernstein Inequality}}:  $|s|_{Y}\leq
C'\,n^{r}\,\|s\|_{X}$, \;for all $s\in X_{n}$,\;$n=1,2,.....$
\\

Our claim now is: Assuming that the modulation space $M_{\infty}$
is the space in which approximation takes place, and using
$\Sigma_{N}$ as the subset of $X$ in which the approximants are
seeked. Then Bernstein inequality holds for our working space.
Before proving our new result let we define $\Sigma_{N}$ using
Gabor frames as follows:
\begin{equation}\label{eq:segma form}
\Sigma_{N}(\mathcal{D})=
    \{s\in M_{\infty};\;s=\sum_{(k,n)\in F}c_{kn}\mathbf{T}_{\alpha k}\mathbf{M}_{\beta n}g,
    \;c_{kn}\in \mathbb{C},\;card(F)\leq N\}
\end{equation}

\begin{prop}\label{prop:Result1-Berns.Modul}
 Let $1\leq p\leq q<\infty$, $\alpha>0$ and $g\in M_{1}$, let
$\mathcal{D}:=
          \{\mathbf{T}_{\alpha k}\mathbf{M}_{\beta n}g;\;k,n\in \mathbb{Z}\}$ be
  a dictionary given by a Gabor frame for $L_{2}(\mathbb{R})$, and
  for $\alpha=(\frac{1}{p}-\frac{1}{q})+1$. Then we have:
  \begin{equation}
  If\;\; s\in \Sigma_{N}(\mathcal{D}),\qquad then\qquad
  \|s\|_{M_{p}}\leq C\,N^{\alpha}\,\|s\|_{M_{q}}\qquad \forall
  N=1,2,....
    \end{equation}
      for some positive constant $C:=C(\alpha,\beta,g)$.
\end{prop}

\begin{proof}
Let $s\in
  \Sigma_{N}$ which has a Gabor expansion:
  \begin{equation}\label{eq:s has Gabor expan}
  s=\sum_{(k,n)\in F}c_{kn}\mathbf{T}_{\alpha k}\mathbf{M}_{\beta
  n}g,
  \end{equation}
  for some indexed set $F$ with $card(F)\leq N$ and some coefficients $c_{kn}$, and for simplicity of notations we let
  $g_{kn}=\mathbf{T}_{\alpha k}\mathbf{M}_{\beta n}g$. Taking the $M_{p}$-norm for $s$ and using
  formula (\ref{eq:M-p,q definition and V_g(f)}) we
have:
  \begin{eqnarray*}
 \|s\|_{M_{p}}^{p}&=&\|\mathcal{V}_{g}s\|_{L_{p}}^{p}\\
                      &\leq&\int_{\mathbb{R}}\int_{\mathbb{R}}\left(\sum_{(k,n)\in F}\big|c_{kn}\mathcal{V}_{g}
                          g_{kn}(x,y)\big|\right)^{p}\;dx\,dy\\
\end{eqnarray*}
and applying the H\"{o}lder's inequality for $\sum_{(k,n)\in
F}\big|c_{kn}\mathcal{V}_{g}g_{kn}(x,y)\big|$, for $1\leq
p<\infty$, and using the fact that $q'\leq p'$ we get
\begin{eqnarray*}
\|s\|_{M_{p}}&\leq& \left(\sum_{(k,n)\in
                        F}|c_{kn}|^{p}\right)^{1/p}.\left(\int_{\mathbb{R}}\int_{\mathbb{R}}\left(\Big[\sum_{(k,n)\in
                   F}\big|\mathcal{V}_{g}g_{kn}(x,y)\big|^{q'}\Big]^{1/q'}\right)^{p}dx\,dy\right)^{1/p}\\
                   &\leq&  N^{\frac{1}{p}}\big(\sup_{(k,n)\in
                          F}|c_{kn}|\,\big)\,.\,N^{\frac{1}{q'}}\left(\int_{\mathbb{R}}
                          \int_{\mathbb{R}}\big(\sup_{(k,n)\in F}\big|\mathcal{V}_{g}g_{kn}(x,y)\big|\big)^{p}dx\,dy
                          \right)^{\frac{1}{p}}\\
                  &=& \,N^{\frac{1}{p}-\frac{1}{q}+1}.\,\big(\sup_{(k,n)\in
                          F}|c_{kn}|\,\big)\,.\|\mathcal{V}_{g}g_{k'n'}\|_{L_{p}},\;\; for\,some\,(k',n')\in F\\
                  &\leq& \,N^{\frac{1}{p}-\frac{1}{q}+1}.\,\big(\sup_{(k,n)\in
                          F}|c_{kn}|\,\big).\,C'\,\|g\|_{M_{p}}
\end{eqnarray*}
where we used part (2) of Theorem (\ref{thm:properties of M-p,q,
Banach and invariant}) in the last inequality. Now, since $g\in
M_{1}$ and $1\leq p<\infty$, then $\|g\|_{M_{p}}\leq
\|g\|_{M_{1}}$ which is finite. Moreover, Since $s\in M_{\infty}$,
Theorem (\ref{thm:Modulation spaces and expansion, and equivalence
relations}) implies that there exist constants $A,\,B>0$ (depends
only on $\alpha,\beta$ and $g$) such that for $p=q=\infty$ we have

\begin{equation*}
\|s\|_{M_{p}}\leq C'.\,
N^{\frac{1}{p}-\frac{1}{q}+1}\,(B.\|s\|_{M_{\infty}})\,.\|g\|_{M_{1}}
\end{equation*}
 Furthermore, we know that $\|s\|_{M_{\infty}}\leq \|s\|_{M_{q}}$ for
$1\leq q<\infty$, hence
\begin{equation*}
\|s\|_{M_{p}}\leq
              \textbf{C}\,N^{\frac{1}{p}-\frac{1}{q}+1}\,\|s\|_{M_{q}}
\end{equation*}
for a constant $\textbf{C}$ independent of $s$, which completes
the proof.

\end{proof}

\textbf{Remark 3.1}.\,(\cite{Devore}). By the monotonicity of the
sequence $(\sigma_{N}(f)_{M_{q}})$, we have the following
equivalence relation:
\begin{equation}\label{eq:equivalence relation series}
  \left(\sum_{N=1}^{\infty}\,\left[N^{\alpha}\,\sigma_{N}
  (f)_{M_{q}}\right]^{\lambda}\frac{1}{N}\right)^{1/\lambda} \asymp
   \left(\sum_{N=0}^{\infty}\,\left[(2^{N})^{\alpha}\,\sigma_{2^{N}}
         (f)_{M_{q}}\right]^{\lambda}\right)^{1/\lambda}
\end{equation}
for each $\alpha>0$ and $0<\lambda<\infty$. Furthermore, this
equivalence means that there exists finite constants
$A_{1},A_{2}>0$, such that
\begin{eqnarray*}
A_{1}\, \left(\sum_{N=1}^{\infty}\,\left[N^{\alpha}\,\sigma_{N}
           (f)_{M_{q}}\right]^{\lambda}\frac{1}{N}\right)^{1/\lambda}
  &\leq&  \left(\sum_{N=0}^{\infty}\,
       \left[(2^{N})^{\alpha}\,\sigma_{2^{N}}(f)_{M_{q}}\right]^{\lambda}\right)^{1/\lambda}\\
  &\leq&\,A_{2}\,  \left(\sum_{N=1}^{\infty}\,\left[N^{\alpha}\,\sigma_{N}
  (f)_{M_{q}}\right]^{\lambda}\frac{1}{N}\right)^{1/\lambda}
\end{eqnarray*}

\begin{thm}\label{thm:newresult-other dirction}
Let $\{\mathbf{T}_{\alpha k}\mathbf{M}_{\beta n}g;\;k,n\in
\mathbb{Z}\}$ be a Gabor frame for $L_{2}(\mathbb{R})$,
$\Sigma_{N}$ as defined in (\ref{eq:segma form}), $1\leq p\leq
q<\infty$, and $\alpha=\frac{1}{p}-\frac{1}{q}+1$. Then we have

  \begin{equation}
   if \qquad
   \sum_{N=1}^{\infty}\left(N^{\alpha}\sigma_{N}(f)_{M_{q}}
      \right)\frac{1}{N}<\infty \qquad then \qquad f\in M_{p}\,.
  \end{equation}
\end{thm}

\begin{proof}
Given a function $f$ belongs to the modulation space $M_{\infty}$
and has a Gabor expansion
\[ f=\sum_{k,n\in \mathbb{Z}}\lambda_{kn}\mathbf{T}_{\alpha k}\mathbf{M}_{\beta n}g
 \]
  for some window function $g\in M_{1}$. Moreover, suppose that the approximation error of
 approximating $f$ by elements from $\Sigma_{N}$ is measured in an
 $M_{q}$-norm and satisfies
\begin{equation}\label{eq:assumption}
\left(\sum_{N=1}^{\infty}\,\left[N^{\alpha
  }\,\sigma_{N}(f)_{M_{q}}\right]\frac{1}{N}\right)<\infty,
 \end{equation}
for any $\alpha>0$, and $1\leq p<\infty$, we need to show that $ \|f\|_{M_{p}}<\infty$.\\
For all $N\in \mathbb{N}$, let $s_{N}\in \Sigma_{2^{N}}$ be a
near-best approximant to $f$ from $\Sigma_{2^{N}}$, i.e.,
\begin{equation}\label{eq:assump2-proof}
\|f-s_{N}\|_{M_{q}}\equiv \sigma_{2^{N}}(f)_{M_{q}}
\end{equation}
Furthermore, we can assume that every $f$ in $M_{\infty}$ can be
written as
\begin{equation*}
     f=\lim_{k\longrightarrow\infty}\sum_{N=1}^{k}(s_{N}-s_{N-1})=\sum_{N=1}^{\infty}(s_{N}-s_{N-1})
\end{equation*}

where $s_{0}=0$. Now, taking the $M_{p}$-norm for $f$ we get
\begin{equation}\label{eq:taking norm-triangle-inq}
    \|f\|_{M_{p}}\leq\sum_{N=1}^{\infty}\|s_{N}-s_{N-1}\|_{M_{p}}
\end{equation}

Since $\Sigma_{2^{N-1}}\subseteq\Sigma_{2^{N}}$. Then
$(s_{N}-s_{N-1})\in \Sigma_{2^{N}}$. So we can apply our result in
Proposition (\ref{prop:Result1-Berns.Modul}) for $(s_{N}-s_{N-1})$
as follows: There exist $C>0$ such that
\begin{equation}
   \|s_{N}-s_{N-1}\|_{M_{p}}\leq C\,(2^{N})^{\alpha}\,\|s_{N}-s_{N-1}\|_{M_{q}}
\end{equation}
for the given $\alpha>0$.

Moreover, from assumption (\ref{eq:assump2-proof}) and from the
monotonicity of $(\sigma_{N}(f)_{M_{q}})$, we have
\begin{equation}\label{eq:S_N - S_n-1 relation}
\|s_{N}-s_{N-1}\|_{M_{q}}\leq C'\,\sigma_{2^{N-1}}(f)_{M_{q}},\;
for\;some\;constant\; C'>0.
\end{equation}
Back to Equation (\ref{eq:taking norm-triangle-inq}), using the
previous inequality, we have,
\begin{equation*}
 \|f\|_{M_{p}}\leq C'_{2}\,\sum_{N=0}^{\infty}\,\left[(2^{N})^{\alpha
                         }\,\sigma_{2^{N}}(f)_{M_{q}}\right]
\end{equation*}
Using Remark (3.1) for $\lambda=1$ and our assumption in
(\ref{eq:assumption}) above we get
\begin{eqnarray*}
    \|f\|_{M_{p}} \leq C''_{2}\left(\sum_{N=1}^{\infty}\,\left[N^{\alpha
                        }\,\sigma_{N}(f)_{M_{q}}\right]\frac{1}{N}\right)<\infty
\end{eqnarray*}
Thus, our proof is complete and we have $ \|f\|_{M_{p}}<\infty$,
i.e., $f\in M_{p}$\,.
\end{proof}


\section{\textbf{Characterization of Modulation Spaces} $M_{p,q}$}
\label{sec:charac.of Modul.spaces M-p,q}
\ \ \ \ Now it is natural according to the later results that one
thinking in the characterization of the modulation spaces
$M_{p,q}$ using the Gabor frames and nonlinear approximation, and
if it is possible to characterize an arbitrary function $f$
belonging to some modulation space like $M_{\infty}$ using
elements from the same space, knowing that its approximation error
measured in the $M_{p,q}$-norm satisfying some condition.

Let $\mathcal{D}:=\{\mathbf{T}_{\alpha k}\mathbf{M}_{\beta
n}g;\,k,n\in \mathbb{Z}\}$ \emph{be  a dictionary given by a Gabor
frame for} $L_{2}(\mathbb{R})$, \emph{and define for all} $N\in
\mathbb{N}$:
\begin{equation}\label{eq:segma general form in X}
       \Sigma_{N}(\mathcal{D})=\{s\in M_{\infty};\, s=\sum_{(k,n)\in F}c_{kn}
       \mathbf{T}_{\alpha k}\mathbf{M}_{\beta n}g,
         \;c_{kn}\in \mathbb{C},\;card(F)\leq N\}
      \end{equation}
  \begin{equation}
   \sigma_{N}(f)_{M_{p,q}}=\inf_{s\in \Sigma_{N}}\|f-s\|_{M_{p,q}}
  \end{equation}

We can solve this problem in one direction. More precisely, for
$\Sigma_{N}$ defined in Equation (\ref{eq:segma general form in
X}), we have

\begin{thm}\label{thm:newresult-other dirction in M-p1,q1}
Let $\{\mathbf{T}_{\alpha k}\mathbf{M}_{\beta n}g;\;k,n\in
\mathbb{Z}\}$ be a Gabor frame for $L_{2}(\mathbb{R})$,
$\Sigma_{N}$ as defined in Equation (\ref{eq:segma general form in
X}), $1\leq p_{1}\leq p<\infty$, $1\leq q_{1}\leq q<\infty$, and
$\alpha=(\frac{1}{p_{1}}+\frac{1}{q_{1}})-(\frac{1}{p}+\frac{1}{q})+2$.
Then,
  \begin{equation}
  if\qquad \sum_{N=1}^{\infty}\left(N^{\alpha}
      \sigma_{N}(f)_{M_{p,q}}\right)\frac{1}{N}<\infty ,  \qquad
      then\;\;\;f\in M_{p_{1},q_{1}}
  \end{equation}
\end{thm}
Before starting the proof of the theorem, we must show that the
Bernstein estimate theorem can be applied to these modulation
spaces under our new definition of $\Sigma_{N}$ and a given rate
$\alpha>0$.

\begin{prop}(\textbf{Bernstein Inequality and Modulation Spaces $M_{p,q}$}).
\label{prop:Result1-Berns.Modul M-pq}
\\  Let $1\leq p_{1}\leq p<\infty$, $1\leq q_{1}\leq q<\infty$ and $\alpha>0$.
Let $g\in M_{1}$, and $\mathcal{D}:=\{\mathbf{T}_{\alpha
k}\mathbf{M}_{\beta n}g;\;k,n\in \mathbb{Z}\}$ be
  a dictionary given by a Gabor frame for $L_{2}(\mathbb{R})$. Then for
  $\alpha=(\frac{1}{p_{1}}+\frac{1}{q_{1}})-(\frac{1}{p}+\frac{1}{q})+2$, we
  have:
  \begin{equation}
  if\;\; s\in \Sigma_{N}(\mathcal{D})\,,\qquad then\qquad
  \|s\|_{M_{p_{1},q_{1}}}\leq C\,N^{\alpha}\,\|s\|_{M_{p,q}},\qquad \forall
  N=1,2,....
    \end{equation}
    for some positive constant $C:=C(\alpha,\beta,g)$.
\end{prop}

\begin{proof}
Let $s\in
  \Sigma_{N}$. Then
  \begin{equation}\label{eq:s has Gabor expan in M-p,q}
  s=\sum_{(k,n)\in F}c_{kn}\mathbf{T}_{\alpha k}\mathbf{M}_{\beta
  n}g=\sum_{k}\sum_{n}c_{kn}\mathbf{T}_{\alpha k}\mathbf{M}_{\beta
  n}g
  \end{equation}
  for some indexed set $F$ with $card(F)\leq N$ and some coefficients
  $c_{kn}$. For simplicity of notations we let $g_{kn}=\mathbf{T}_{\alpha k}\mathbf{M}_{\beta n}g$.
  Now, taking the $M_{p_{1},q_{1}}$-norm for $s$ and using formula (\ref{eq:M-p,q definition and V_g(f)})
  we get
\begin{eqnarray*}
\|s\|_{M_{p_{1},q_{1}}}^{q_{1}}&=& \|\mathcal{V}_{g}s\|_{L_{p_{1},q_{1}}}^{q_{1}}\\
                     &\leq& \int_{\mathbb{R}}\left(\int_{\mathbb{R}}\left[\sum_{k}\sum_{n}\big|c_{kn}\big|\big|\mathcal{V}_{g}
                            g_{kn}(x,y)\big|\right]^{p_{1}}\;dx\right)^{q_{1}/p_{1}}\,dy
\end{eqnarray*}
Now, keeping in mind that we are working over a finite indexed set
$F$, using the H\"{o}lder's inequality over the internal index,
and $\frac{1}{p_{1}}+\frac{1}{p'_{1}}=1$, we have for each $x,y$
\[ \sum_{k}\left(\sum_{n}\big|c_{kn}\big|\big|\mathcal{V}_{g}
                            g_{kn}(x,y)\big|\right) \leq\sum_{k}\left( \Big(\sum_{n}\big|c_{kn}\big|^{p_{1}}
                            \Big)^{\frac{1}{p_{1}}} .\,\Big(\sum_{n}\big|\mathcal{V}_{g}
                            g_{kn}(x,y)\big|^{p'_{1}}\Big)^{\frac{1}{p'_{1}}} \right)
\]

Again, using the H\"{o}lder's inequality over the external index,
and $\frac{1}{q_{1}}+\frac{1}{q'_{1}}=1$, the right hand side of
the last inequality will be less than or equal to
\begin{equation}\label{eq:Holders ineq in second time}
\left(\sum_{k}\left[\sum_{n}\big|c_{kn}\big|^{p_{1}}
                                 \right]^{\frac{q_{1}}{p_{1}}}\right)^{\frac{1}{q_{1}}} .\,
                                 \left(\sum_{k}\left[\sum_{n}\big|\mathcal{V}_{g}
                              g_{kn}(x,y)\big|^{p'_{1}}\right]^{\frac{q'_{1}}{p'_{1}}} \right)^{\frac{1}{q'_{1}}}
\end{equation}
Now, simplifying the first part of the expression in
(\ref{eq:Holders ineq in second time}), we get
\begin{equation*}
 \left(\sum_{k}\left[\sum_{n}\big|c_{kn}\big|^{p_{1}}
              \right]^{\frac{q_{1}}{p_{1}}}\right)^{\frac{1}{q_{1}}}
       \leq N^{\frac{1}{p_{1}}+\frac{1}{q_{1}}}.\,\Big(\sup_{(k,n)\in
             F}\big|c_{kn}\big|\,\Big)
\end{equation*}

But since $p_{1}\leq p$, then $p'\leq p'_{1}$, so
   \[\Big[\sum_{n}\big|\mathcal{V}_{g}g_{kn}(x,y)\big|^{p'_{1}}\Big]^{1/p'_{1}} \leq
       \Big[\sum_{n}\big|\mathcal{V}_{g}g_{kn}(x,y)\big|^{p'}\Big]^{1/p'}
   \]

and we get
\begin{eqnarray*}
\left(\sum_{k}\left[\sum_{n}\big|\mathcal{V}_{g}
     g_{kn}(x,y)\big|^{p'_{1}}\right]^{\frac{q'_{1}}{p'_{1}}}
                   \right)^{\frac{1}{q'_{1}}}&\leq& \left(\sum_{k}\left[N^{\frac{1}{p'}}\,
                   \sup_{n}\big|\mathcal{V}_{g} g_{kn}(x,y)\big|\right]^{q'_{1}}
                   \right)^{\frac{1}{q'_{1}}}\\
            &\leq& N^{\frac{1}{p'}}\,\left(\sum_{k}\big|\mathcal{V}_{g}
                  g_{kn'}(x,y)\big|^{q'_{1}}
                   \right)^{\frac{1}{q'_{1}}}
\end{eqnarray*}
Since $q_{1}\leq q$, we must have $q'\leq q'_{1}$. Hence we get
\begin{eqnarray*}
\left(\sum_{k}\big|\mathcal{V}_{g}g_{kn'}(x,y)\big|^{q'_{1}}\right)^{1/q'_{1}}
            &\leq& \left(\sum_{k}\big|\mathcal{V}_{g}g_{kn'}(x,y)\big|^{q'}\right)^{1/q'}\\
             &=& N^{\frac{1}{q'}}\big|\mathcal{V}_{g}
                  g_{k'n'}(x,y)\big|
\end{eqnarray*}
Therefore
\begin{equation*}
\left(\sum_{k}\left[\sum_{n}\big|\mathcal{V}_{g}
                  g_{kn}(x,y)\big|^{p'_{1}}\right]^{\frac{q'_{1}}{p'_{1}}}
                   \right)^{\frac{1}{q'_{1}}}
            \leq N^{2-(\frac{1}{p}+\frac{1}{q})}\,\big|\mathcal{V}_{g}
                  g_{k'n'}(x,y)\big|
\end{equation*}
Combining our results we get
\begin{eqnarray*}
\|s\|_{M_{p_{1},q_{1}}} &\leq&
\left(\int_{\mathbb{R}}\left(\int_{\mathbb{R}}\left[\sum_{k}\sum_{n}\big|c_{kn}\big|\big|\mathcal{V}_{g}
                                        g_{kn}(x,y)\big|\right]^{p_{1}}\;dx\right)^{q_{1}/p_{1}}\,dy\right)^{\frac{1}{q_{1}}}\\
                                 & &\qquad\qquad \left(\sum_{k}\left[\sum_{n}\big|\mathcal{V}_{g}
                                        g_{kn}(x,y)\big|^{p'_{1}}\right]^{\frac{q'_{1}}{p'_{1}}} \right)^{\frac{1}{q'_{1}}}
                                     \Big]^{p_{1}}\;dx\Bigg)^{q_{1}/p_{1}}\,dy\Big)^{\frac{1}{q_{1}}}\\
                              &\leq& N^{(\frac{1}{p_{1}}+\frac{1}{q_{1}})-(\frac{1}{p}+\frac{1}{q})+2}.\,\sup_{(k,n)\in
                                          F}\big|c_{kn}\big|.\,\left(\int_{\mathbb{R}}\left(\int_{\mathbb{R}}
                                          \big|\mathcal{V}_{g}g_{k'n'}(x,y)\big|^{p_{1}}dx\right)^{q_{1}/p_{1}}dy\right)^{1/q_{1}}
\end{eqnarray*}
Using Equation (\ref{eq:M-p,q definition and V_g(f)}) and Theorem
(\ref{thm:Modulation spaces and expansion, and equivalence
relations}) for $p,q< \infty$ we have
\begin{eqnarray*}
\|s\|_{M_{p_{1},q_{1}}}
                   &\leq& N^{(\frac{1}{p_{1}}+\frac{1}{q_{1}})-(\frac{1}{p}+\frac{1}{q})+2}.\,
                           \|s\|_{M_{\infty}}.\,\|\mathcal{V}_{g}g_{k'n'}\|_{L_{p_{1},q_{1}}}\\
                   &\leq& C\,N^{(\frac{1}{p_{1}}+\frac{1}{q_{1}})-(\frac{1}{p}+\frac{1}{q})+2}.\,
                          \|s\|_{M_{p,q}}.\,\|g_{k'n'}\|_{M_{p_{1},q_{1}}}
\end{eqnarray*}

for some positive constant $C:=C(\alpha,\beta,g)$. Moreover, using
part (2) of Theorem (\ref{thm:properties of M-p,q, Banach and
invariant}) and our assumption that $g\in M_{1}$ and $1\leq
p_{1},q_{1}<\infty$ we get our result.
\end{proof}

\textbf{\\Proof of Theorem (\ref{thm:newresult-other dirction in
M-p1,q1})}. \\Given a function $f$ belongs to the modulation space
$M_{\infty}$ and has a Gabor expansion
\[ f=\sum_{k,n\in \mathbb{Z}}\lambda_{kn}\mathbf{T}_{\alpha k}\mathbf{M}_{\beta n}g
 \]
  for some window function $g\in M_{1}$.
  Suppose that the approximation error of
 approximating $f$ by elements from $\Sigma_{N}$ is measured in an
 $M_{p,q}$-norm and satisfies
 \begin{equation}\label{eq:assumption in M-p,g}
\left(\sum_{N=1}^{\infty}\,\left[N^{\alpha
  }\,\sigma_{N}(f)_{M_{p,q}}\right]\frac{1}{N}\right)<\infty
\end{equation}
we need to show that $ \|f\|_{M_{p_{1},q_{1}}}<\infty$.\\
Let $s_{N}\in \Sigma_{2^{N}}$ be a near-best approximant to $f$
from $\Sigma_{2^{N}}$, for all $N\in \mathbb{N}$, i.e.,
\begin{equation}\label{eq:assump2-proof in M-p1,q1 thm}
\|f-s_{N}\|_{M_{p,q}}\equiv \sigma_{2^{N}}(f)_{M_{p,q}}
\end{equation}
Since every $f$ in $M_{\infty}$ can be written as
\begin{equation*}
     f=\sum_{N=1}^{\infty}(s_{N}-s_{N-1})
\end{equation*}

where $s_{0}=0$. Taking the $M_{p_{1},q_{1}}$-norm for $f$ we get
\begin{equation}\label{eq:taking norm-triangle-inq in M-p1,q1 thm}
    \|f\|_{M_{p_{1},q_{1}}}\leq\sum_{N=1}^{\infty}\|s_{N}-s_{N-1}\|_{M_{p_{1},q_{1}}}
\end{equation}

Since $(s_{N}-s_{N-1})\in \Sigma_{2^{N}}$. So, from Proposition
(\ref{prop:Result1-Berns.Modul M-pq}), there exist $C>0$ such that
\begin{equation}
   \|s_{N}-s_{N-1}\|_{M_{p_{1},q_{1}}}\leq C\,(2^{N})^{\alpha}\,\|s_{N}-s_{N-1}\|_{M_{p,q}}
\end{equation}
for the given $\alpha$. Moreover, from assumption
(\ref{eq:assump2-proof in M-p1,q1 thm}) and from the monotonicity
of $(\sigma_{N}(f)_{M_{p,q}})$, we have
\begin{equation}\label{eq:S_N - S_n-1 relation in M-p1,q1 thm}
\|s_{N}-s_{N-1}\|_{M_{p,q}}\leq
C'\,\sigma_{2^{N-1}}(f)_{M_{p,q}},\; for\;some\;constant\; C'>0.
\end{equation}

Back to Equation (\ref{eq:taking norm-triangle-inq in M-p1,q1
thm}), and using the previous inequality we have,
\begin{equation*}
 \|f\|_{M_{p_{1},q_{1}}}\leq C'_{2}\,\sum_{N=0}^{\infty}\,\left[(2^{N})^{\alpha
                         }\,\sigma_{2^{N}}(f)_{M_{p,q}}\right]
\end{equation*}
and using Remark (3.1) for $\lambda=1$ and our assumption in
(\ref{eq:assumption in M-p,g}) above we get
\begin{eqnarray*}
    \|f\|_{M_{p_{1},q_{1}}} \leq C''_{2}\left(\sum_{N=1}^{\infty}\,\left[N^{\alpha
                        }\,\sigma_{N}(f)_{M_{p,q}}\right]\frac{1}{N}\right)<\infty
\end{eqnarray*}
Thus, $ \|f\|_{M_{p_{1},q_{1}}}<\infty$, and $f\in
M_{p_{1},q_{1}}$. $\hspace{2.5in}\quad\square$



\end{document}